\documentclass[11pt]{article}

\usepackage{amssymb}
\usepackage{amsmath}

\usepackage{graphicx}
\graphicspath{%
    {converted_graphics/}% inserted by PCTeX
    {/}% inserted by PCTeX
    {Excptwmf/}% inserted by PCTeX
}

\renewcommand{\Re}{\mathbb{R}}

\begin{document}

\title{A Counterexample in Ito Integration Theory}
\author{Lars Tyge Nielsen \\ Department of Mathematics \\ Columbia University}
\date{\today}
\maketitle

\begin{abstract}
Ito's Lemma implies that if $W$ is a Wiener process and $f$ is a twice continuously differentiable function, then
the process $f(W)$ is the sum of a time integral and an Ito integral. The Ito integrand is not necessarily locally square integrable.
This note provides a counterexample.
\end{abstract}

\section{Introduction}

Occasionally, an exposition of the Ito integral, or the stochastic integral with respect to a Wiener process, will limit the integrands to locally square integrable processes rather than the more general case of \emph{pathwise} locally square integrable processes. See  
Bingham and Kiesel \cite[1998]{Bingham-Kiesel:98},
Bj{\"o}rk \cite[2009]{Bjoerk:09},
Etheridge \cite[2002]{Etheridge:02},
Evans \cite[2013]{Evans:13},
Mikosch \cite[1998]{Mikosch:98},
and Shreve \cite[2004]{Shreve:04}. 
This is done in part for simplicity, and in part to make sure the integrals are martingales. 
Noble aims, but
unfortunately, Ito's Lemma does not hold if the integral is only defined for locally square integrable integrands.

Ito's Lemma implies that if $W$ is a Wiener process and $f$ is a twice continuously differentiable function, then
the process $f(W)$ is the sum of a time integral and an Ito integral. 
However, the Ito integrand is not necessarily locally square integrable.
This note provides a counterexample.

\section{Background}

A \emph{setup} is a quintuple $(\Omega, \mathcal{F},P,F,W)$ consisting of 
\begin{itemize}
\item
a complete probability space $(\Omega, \mathcal{F},P)$,
\item
an augmented filtration $F = (\mathcal{F}_{t})_{t \in [0,\infty)}$, and 
\item
a standard Wiener process $W$ with respect to $F$. 
\end{itemize}
We shall assume that a specific setup has been chosen. 

All processes in this paper are understood to be one-dimensional.

Let $\mathcal{L}^{1}$ be the set of progressively measurable (or at the very least, measurable and adapted) processes $a$ that are \emph{pathwise locally integrable}, which means that for every $t \in [0,\infty)$,
\[ \int_{0}^{t} | a | \, ds  < \infty \]
with probability one.
Let
$\mathcal{L}^{2} \subset \mathcal{L}^{1}$ denote the subset of processes $b$ that are \emph{pathwise locally square integrable}, which means that for every $t \in [0,\infty)$,
\[ \int_{0}^{t} b^{2} \, ds  < \infty \]
with probability one.
Finally, let
$\mathcal{H}^{2} \subset \mathcal{L}^{2}$ denote the subset of processes $b \in \mathcal{L}^{2}$ that are \emph{locally square integrable}, meaning that for every $t \in [0,\infty)$,
\[ E \int_{0}^{t} b^{2} \, ds < \infty \]
If $b \in \mathcal{H}^{2}$, then the Ito integral process
\[ \int_{0}^{t}b \, dW \]
is a square integrable martingale.

An \emph{Ito process} is a process $X$ of the form,
\[ X(t) = X(0) + \int_{0}^{t} a \, ds + \int_{0}^{t} b \, dW \]
where $a \in \mathcal{L}^{2}$ and $b \in \mathcal{L}^{2}$.

Given a twice continuously differentiable function
$f: \Re \rightarrow \Re$, Ito's lemma says that (1) $f(X)$ is an Ito process, and (2) the relevant integrand processes are
$f'(X)a + (1/2)
f''(X) b^{2}$ and $f'(X) b$.
Written in integral form,
\[ f(X(t)) = f(X(0)) +  \int_{0}^{t} \left(f'(X)a + \frac{1}{2} f''(X) b^{2} \right) ds + \int_{0}^{t} f'(X) b \, dW \]
The problem is that when the Ito integral is defined only for integrand processes in $\mathcal{H}^{2}$, the Ito integral part of $f(X(t))$ may well be undefined, and thus Ito's Lemma fails.

In the simple case where $a = 0$ and $b = 1$, the corresponding Ito process is $X = W$, and the Ito-integral part of $f(X) = f(W)$ is
\[ \int_{0}^{t} f'(W) \, dW \]
We shall construct an example of a function $f$ such that the integrand process $f'(W)$ is not in $\mathcal{H}^{2}$ and, consequently, the integral is not 
well defined when you restrict Ito integration to locally square integrable integrand processes nly.

\section{The Counterexample}

The idea is to find a function $f$ such that $f'$ transforms a standard normal distribution into a Student $t$ distribution with two degrees of freedom, which has infinite variance. It will then transform any mean-zero normal distribution with variance greater than one into a distribution that has infinite variance as well. Hence,
for $t \geq 1$,
\[ E \int_{0}^{t} f'(W(s))^{2} \, ds = \int_{0}^{t} E f'(W(s))^{2} \geq \int_{1}^{t} E f'(W(s))^{2}= \infty \]
implying that the process $f'(W)$ is not in $\mathcal{H}^{2}$.

Let $F$ be the cumulative distribution function of a Student $t$ distribution with two degrees of freedom.
This distribution has density function
\[ d(x) = \frac{1}{(2+x^{2})^{3/2}} \]
and cumulative distribution function is
\[ F(x) = \frac{1}{2}+\frac{x}{2\sqrt{2+x^{2}}} \]
The density function $d$ is strictly positive, symmetric around zero, and continuously differentiable.
The distribution has mean zero and infinite variance. The cumulative distribution function $F$ is continuous and strictly increasing and is a bijection from $\Re$ to the interval $(0,1)$.

Let $N$ be the cumulative distribution function of a standard normal distribution.
This distribution has density function
\[ \phi(x) = \frac{1}{\sqrt{2\pi}} \exp\left(-\frac{x^{2}}{2}\right) \]
The density function $\phi$ is strictly positive, symmetric around zero, and continuously differentiable.
The distribution has mean zero and variance one. The cumulative distribution function $N$ is continuous and strictly increasing and is a bijection from $\Re$ to the interval $(0,1)$.

Now define $h = F^{-1} \circ N: \Re \rightarrow \Re$. This function is a strictly increasing and continuously differentiable bijection with strictly positive derivative.

Observe that $h^{-1} = N^{-1} \circ F$, 
\[ h(0) = F^{-1}(N(0)) = F^{-1}(1/2) = 0 \]
and $h^{-1}(0) = 0$. 
This implies that
$h^{-1}(x) > 0$ when $x > 0$ and
$h^{-1}(x) < 0$ when $x < 0$.

Let $f$ be an indefinite integral of $h$, that is, a differentiable function $\Re \rightarrow \Re$ with $f' = h$. 
Since $h$ is continuously differentiable, $f$ is twice continuously differentiable.

The random variable
$h(W(1))$ has cumulative distribution function $F$, and hence it
follows the Student $t$ distribution with two degrees of freedom.
To see this, observe that $h^{-1} = N^{-1} \circ F$, $F = N \circ h^{-1}$, and
\[ P(h(W(1)) \leq y) = P(W(1) \leq h^{-1}(y)) = N(h^{-1}(y)) = F(y)\]
In particular, $Eh(W(1))^{2} = \infty$.

Next, let's show that $Eh(W(s))^{2} = \infty$ when $s \geq 1$.

For $\sigma > 0$, let $G_{\sigma}$ be the cumulative distribution function of $h(\sigma W(1))^{2}$. In particular,
$G_{1}$ is the cumulative distribution function of $h(W(1))^{2}$.

For $y > 0$,
\begin{align*}
G_{\sigma}(y) 
& = 
P(h(\sigma W(1))^{2} \leq y) \\
& = 
P(-\sqrt{y} \leq h(\sigma W(1)) \leq \sqrt{y}) \\
& = 
P(h^{-1}(-\sqrt{y}) \leq \sigma W(1) \leq h^{-1}(\sqrt{y})) \\
& = 
P(h^{-1}(-\sqrt{y})/\sigma \leq W(1) \leq h^{-1}(\sqrt{y}))/\sigma \\
& = 
N(h^{-1}(\sqrt{y}))/\sigma)-N(h^{-1}(-\sqrt{y})/\sigma)
\end{align*}
Since $h^{-1}(\sqrt{y})) > 0$ and $h^{-1}(-\sqrt{y})) < 0$,
as $\sigma$ increases,
$h^{-1}(\sqrt{y}))/\sigma$ decreases towards zero and
$h^{-1}(-\sqrt{y}))/\sigma$ increases towards zero. 
Hence, as $\sigma$ increases,
$G_{\sigma}(y)$ decreases towards zero and
\[ P(h(\sigma W(1))^{2} > y) = 1-G_{\sigma}(y) \]
increases towards one. 

For a non-negative random variable $X$,
\[ EX = \int_{0}^{\infty} P(X>t) \, dt = \int_{0}^{\infty} P(X \geq t) \, dt \]
See Billingsley \cite[1986, page 282]{Billingsley:86}.

Hence,
\begin{align*}
\int_{0}^{\infty} [1-G_{1}(y)]\, dy
& =
\int_{0}^{\infty} P(h( W(1))^{2} > y)\, dy \\
& = 
Eh( W(1))^{2} \\
& = 
\infty
\end{align*}
and for $\sigma \geq 1$,
\begin{align*}
Eh(\sigma W(1))^{2} 
& = 
\int_{0}^{\infty} P(h(\sigma W(1))^{2} > y)\, dy \\
& = 
\int_{0}^{\infty} [1-G_{\sigma}(y)]\, dy \\
& \geq
\int_{0}^{\infty} [1-G_{1}(y)]\, dy \\
& =
\infty
\end{align*}
For $s > 0$, $W(s)$ and $\sqrt{s} \, W(1)$ are both normal with mean zero and variance $s$.
Hence, for $t > 1$,
\[ E \int_{1}^{t} h(W(s))^{2} \, ds = 
\int_{1}^{t} E h(\sqrt{s} \, W(1))^{2} \, ds = \infty \] 
and
\begin{align*}
E \int_{0}^{t} h(W(s))^{2} \, ds
& =
\int_{0}^{t} Eh(W(s))^{2} \, ds \\
& =
\int_{0}^{1} Eh(W(s))^{2} \, ds + \int_{1}^{t} Eh(W(s))^{2} \, ds \\
& =
\int_{0}^{1} Eh(W(s))^{2} \, ds + \infty \\
& =
\infty
\end{align*}
This implies that $f'(W) = h(W)$ does not belong to $\mathcal{H}^{2}$, and consequently,
the integral 
\[ \int_{0}^{t} f'(W) \, dW \]
is not well defined, if you choose to define the Ito integral only for locally square integrable integrands.

%\bibliographystyle{plain}
%\bibliography{bookrefs,bookrefsadd,consolidated}

\end{document}